\documentclass[11pt,reqno]{amsart}
\usepackage{amssymb,amscd,amsbsy,mathrsfs}
\newcommand{\lb}{\linebreak}

\newcommand{\e}{\varepsilon}

\newcommand{\z}{\zeta}

\newcommand{\s}{\sigma}
\renewcommand{\t}{\tau}

\newcommand{\f}{\varphi}

\newcommand{\D}{\Delta}

\newcommand{\h}{{\mathscr H}}

\newcommand{\X}{{\mathscr X}}
\newcommand{\Y}{{\mathscr Y}}

\newcommand{\T}{{\Bbb T}}

\newcommand{\R}{{\Bbb R}}

\newcommand{\0}{{\boldsymbol{0}}}

\newcommand{\bs}{\boldsymbol}

\newcommand{\bS}{{\boldsymbol S}}

\newcommand{\rf}[1]{(\ref{#1})}

\newcommand{\df}{\stackrel{\mathrm{def}}{=}}

\newcommand{\supp}{\operatorname{supp}}

\newcommand{\trace}{\operatorname{trace}}

\newcommand{\const}{\operatorname{const}}

\newcommand{\eeq}{\end{equation}}
\newcommand{\beq}{\begin{equation}}
\newcommand{\bay}{\begin{eqnarray}}
\newcommand{\ba}{\begin{align*}}
\newcommand{\ea}{\end{align*}}
\newcommand{\ey}{\end{eqnarray}}
\newcommand{\bey}{\begin{eqnarray*}}
\newcommand{\eey}{\end{eqnarray*}}

\newcommand{\be}{\infty}

\newcommand{\bl}{\blacksquare}

\newcommand{\Pf}{{\bf Proof. }}

\newtheorem{thm}{\hspace{\parindent}Theorem}[section]

\newtheorem{lem}[thm]{\hspace{\parindent}Lemma}

\usepackage{amsmath}

\theoremstyle{remark}

\newtheorem*{rem*}{Remark}

\newcommand\fM{\frak M}

\newcommand\dg{\frak D}

\newcommand\mB{\mathcal{B}}

\newcommand{\OL}{{\rm OL}}

\begin{document}

\

\newcommand{\vse}{\vspace{.2in}}
\numberwithin{equation}{section}

\title{Krein's trace formula for unitary operators and operator Lipschitz functions}
\author{A.B. Aleksandrov and V.V. Peller}
\thanks{The first author is partially supported by the RFBR grant 14-01-00198; the second author is partially supported 
by the NSF grant DMS 1300924.}

\begin{abstract}
The main result of this paper is a description of the space of functions on the unit circle, for which Krein's trace formula
holds for arbitrary pairs of unitary operators with trace class difference. This space coincides with the space of operator Lipschitz
functions on the unit circle.
\end{abstract}

\maketitle

\tableofcontents

\

\setcounter{section}{0}
\section{\bf Introduction}
\setcounter{equation}{0}
\label{Vved}

\

The {\it spectral shift function} for  pairs of selfadjoint
operators was introduced in the paper by I.M. Lifshits  \cite{Li}.
In the same paper a trace formula for the difference of functions of the perturbed operator and the unperturbed operator was established. Ideas by Lifshits were developed in the paper by M.G. Krein  \cite{Kr1}, in which the spectral shift function $\bs{\xi}$ in $L^1(\R)$
was defined for arbitrary pairs of self-adjoint operators $A$ and $B$ with $A-B$ in trace class,
and the Lifshits trace formula 
\bay
\label{fsLK}
\trace\big(f(A)-f(B)\big)=\int_\R f'(t)\bs{\xi}(t)\,dt,
\ey
was proved in the considerably more general situation when the derivative of $f$ is the Fourier transform of
a complex Borel measure on $\R$.

Later in \cite{Pe2} and \cite{Pe4} trace formula  \rf{fsLK} was extended to arbitrary functions $f$
in the Besov class $B_{\be,1}^1(\R)$ (see \cite{Pee} for definitions of the Besov classes).

On the other hand, it is obvious that the right-hand side of  \rf{fsLK} is well defined for an arbitrary Lipschitz function $f$.
M.G. Krein posed in \cite{Kr1} the question of whether it is possible to generalize formula  \rf{fsLK} to
the case of arbitrary Lipschitz functions. It turned out that the answer to this question is negative:
Yu.B. Farforovskaya in  \cite{F2} constructed an example of a Lipschiz function $f$ and 
self-adjoint operators $A$ and $B$ such that $A-B$ belongs to trace class  $\bS_1$
but $f(A)-f(B)\notin\bS_1$.

Thus, the question of applicability of trace formula  \rf{fsLK} actually splits in two
separate questions:

\medskip

(a) {\it For what functions $f$ on $\R$ does the implication
$$
A-B\in\bS_1\quad\Longrightarrow\quad f(A)-f(B)\in\bS_1
$$
hold for not necessarily bounded self-adjoint operators  $A$ and $B$?}

(b) {\it If $f$ satisfies condition} (a){\it, is it true  that the left-hand side and the right-hand side of equality 
{\em\rf{fsLK}} coincide?}

\medskip

It is well known (see, for example, recent survey \cite{AP}, Th. 3.6.5) that a function  $f$ on $\R$
satisfies (a) if and only if it is {\it operator Lipschitz}, i.e., the inequality
$$
\|f(A)-f(B)\|\le\const\|A-B\|
$$
holds for arbitrary (bounded or unbounded) self-adjoint operators $A$ and $B$.

The fact that a Lipschitz function does not have to be operator Lipschitz was established in the paper by Yu.L. Farforovskaya \cite{F1}. Then it was proved in \cite{JW} that an operator Lipschitz function is differentiable everywhere. 
This readily implies the earlier results of papers  \cite{Mc} and \cite{Ka}: the function $x\mapsto|x|$
is not operator Lipschitz. Note also that an operator Lipschitz function is not necessarily 
continuously differentiable which was proved in \cite{KS}.
 Necessary conditions for operator Lipschitzness were obtained in  \cite{Pe2} and \cite{Pe4}.
 These necessary conditions are based on the description  \cite{Pe1} of trace class Hankel operators
 (see also \cite{Pe3}).

 We refer the reader to the survey article \cite{AP}, which provides a detailed analysis of sufficient conditions and necessary conditions
 for operator Lipschitzness.
  
A positive answer to question (b) was obtained in a recent paper \cite{Pe5}: formula \rf{fsLK}
is true for arbitrary operator Lipschitz fubction  $f$. Thus, the class of functions, for which trace formula 
 \rf{fsLK} holds for all self-adjoint operators $A$ and $B$ with trace class difference 
coincides with the class of operator Lipschitz functions.

In this paper we propose a solution to a similar problem for functions of unitary operators.

Spectral shift function for pairs of unitary operators with trace class difference was introduced
in M.G. Krein's paper  \cite{Kr2} (see also \cite{Kr3}, in which a detailed presentation of the results is given).
Let $U$ and $V$ be unitary operators  with $U-V$ in trace class. Then there exists an integrable function
$\bs{\xi}$ on the unit circle $\T$ (called {\it a spectral shift function} for the pair $(U,V)\:$)
such that the trace formula
\bay
\label{LKun}
\trace\big(f(U)-f(V)\big)=\int_\T f'(\z)\bs{\xi}(\z)\,d\z
\ey
holds for sufficiently nice functions $f$.
 In contrast with the case of self-adjoint operators, the function $\bs{\xi}$ is not determined uniquely
 by the pair $(U,V)$; it is determined uniquely up to a constant additive. Thus, it is reasonable to require that
 the mean value of $\bs{\xi}$ over $\T$ to be $0$.

In \cite{Kr2} M.G. Krein  showed that trace formula \rf{LKun} holds
if the derivative $f'$ has absolutely convergent Fourier series. In  \cite{Pe2} the trace formula
was extended to the functions $f$ in the Besov space $B_{\be,1}^1(\T)$.

Note that as in the case of self-adjoint operators, a function $f$ on the unit circle $\T$ takes trace class perturbations to trace class increments, i.e.,
$$
U-V\in\bS_1\quad\Longrightarrow\quad f(U)-f(V)\in\bS_1
$$
if and only if  $f$ is an {\it operator Lipschitz function}, i.e.,
$$
\|f(U)-f(V)\|\le\const\|U-V\|
$$
for all unitary operators $U$ and $V$. In the space $\OL(\T)$ of operator Lipschitz functions we introduce  
the natural seminorm 
$$
\|f\|_{\OL}\df\sup\frac{\|f(U)-f(V)\|}{\|U-V\|},
$$
where the supremum is taken over all unitary operators $U$ and $V$ such that $U\ne V$.

The principal result of this paper is obtained in \S\:\ref{Unitar}; it says that  trace formula \rf{LKun} 
holds for an arbitrary operator Lipschitz function $f$. Clearly, this is the maximal class of functions with this property.
This result will imply the following amusing fact: the function
$$
\z\mapsto\trace\big(f(\z U)-f(\z V)\big)
$$
is continuous on $\T$ for an arbitrary operator Lipschitz function $f$ and an arbitrary pair $(U,V)$
of unitary operators with trace class difference.

Note that the proof obtained in  \cite{Pe5} for functions of self-adjoint operators does not extend
to the case of unitary operators because it is based on a result of the paper \cite{KPSS} on the 
differentiability of operator functions in the Hilbert--Schmidt norm. We do not know whether an analogue 
of their result holds in the case of functions of unitary operators.

Instead, we use in this paper the differentiability of the corresponding operator functions
in the strong operator topology which will be established in \S\:\ref{Difsil}.

In  \S\:\ref{Dvoiopi} we give a brief introduction in double operator integrals. Moreover, we obtain a general
trace formula, which we are going to apply in \S\:\ref{Unitar} for the proof of the main result.

Finally, in \S\:\ref{Samosop} we consider briefly an alternative approach in the case of self-adjoint operators,
which in contrast with the proof obtained in \cite{Pe5} uses the differentiability of operator functions in the strong
operator topology rather than in the Hilbrt--Schmidt norm.

\

\section{\bf Double operator integrals and Schur multipliers}
\setcounter{equation}{0}
\label{Dvoiopi}

\

Double operator integrals appeared in the paper by Yu.L. Daletskii and S.G. Krein  \cite{DK}. Then 
Birman and Solomyak in \cite{BS1}, \cite{BS2} and \cite{BS4} developed a beautiful theory of
double operator integrals.

Let  $(\X,E_1)$ and $(\Y,E_2)$ be spaces with spectral measures $E_1$ and $E_2$
on a Hilbert space  $\h$, and let  $\Phi$ be a bounded measurable function on $\X\times\Y$.
Double operator integrals are expressions of the form
\bay
\label{doi}
\int\limits_\X\int\limits_\Y\Phi(x,y)\,d E_1(x)T\,dE_2(y).
\ey

A starting point in the papers by Birman and Solomyak was the case when $T$ is a Hilbert--Schmidt operator;
in this case double operator integrals can be defined for arbitrary bounded measurable functions $\Phi$. 

We are not going to consider here the case of Hilbert--Schmidt operators and refer the reader to the survey article
\cite{AP}, Chapter II, where double operator integrals are studied in details.

To define double operator integrals of the form \rf{doi} for arbitrary bounded operators $T$,
we have to impose restrictions on the function $\Phi$. In fact, double operator infegrals
can be defined for all bounded operators $T$ under the assumption that $\Phi$ belongs to the class of
 {\it Schur multipliers} $\fM(E_1,E_2)$ with respect to the spectral measures $E_1$ and $E_2$.
 The space $\fM(E_1,E_2)$ admits different descriptions, see  \cite{Pe2}, \cite{Pi} and \cite{AP}. 

We give here one of the descriptions: {\it $\Phi\in\fM(E_1,E_2)$ if and only if $\Phi$ belongs to
the Haagerup tensor product $L^\be(E_1)\otimes_{\rm h}L^\be(E_2)$, i.e., $\Phi$ admits
a representation
\bay
\label{predstava}
\Phi(x,y)=\sum_n\f_n(x)\psi_n(y),
\ey
where $\f_n$ and $\psi_n$ are measurable functions satisfying the condition}
$$
\Big\|\sum_n|\f_n|^2\Big\|_{L^\be(E_1)}\le\const\|\Phi\|_{\fM(E_1,E_2)}
\quad\mbox{and}\quad
\Big\|\sum_n|\psi_n|^2\Big\|_{L^\be(E_2)}\le\const\|\Phi\|_{\fM(E_1,E_2)},
$$
where $\|\Phi\|_{\fM(E_1,E_2)}$ is the norm of the transformer
$$
T\mapsto\iint\Phi\,dE_1T\,dE_2
$$
on the space of operators on Hilbert space.
Then the following equality holds
\bay
\label{doitproh}
\int\limits_\X\int\limits_\Y\Phi(x,y)\,d E_1(x)T\,dE_2(y)=
\sum_n\Big(\int\f_n\,dE_1\Big)T\Big(\int\psi_n\,dE_2\Big),
\ey
where the series on the right-hand side of the equality converges 
in the weak operator topology, and its sum does not depend on
a representation of the form \rf{predstava}.

If $\Phi\in\fM(E_1,E_2)$ and $T$ is a trace class operator, then
the double operator integral \rf{doi} also belongs to trace class and
the inequality
\bay
\label{yadotse}
\left\|\,\int\limits_\X\int\limits_\Y\Phi(x,y)\,d E_1(x)T\,dE_2(y)\right\|_{\bS_1}
\le\|\Phi\|_{\fM(E_1,E_2)}\|T\|_{\bS_1}
\ey
holds.

Suppose that $f$ is an operator Lipschitz function on the unit circle $\T$.
Consider its divided difference on $\T\times\T$:
$$
(\dg f)(\z,\t)\df\left\{\begin{array}{ll}\frac{f(\z)-f(\t)}{\z-\t},&\z\ne\t,\\[.2cm]
f'(\z),&\z=\t
\end{array}\right.
$$
(by virtue of results of \cite{JW}, an operator Lipschitz function on the circle is differentiable at each point).
It is well known that in this case the divided difference $\dg f$ is a Schur multiplier for every
Borel spectral measures $E_1$ and $E_2$. In fact, the converse is also true: if a function $f$ on $\T$
is differentiable everywhere and $\dg f$ is a Schur multiplier for all Borel spectral measures, then
$f$ is operator Lipschitz (see, for example, survey \cite{AP}, Th. 3.3.6). Moreover, the following
equality holds:
$$
\|f\|_{\OL}=\sup\|\dg f\|_{\fM(E_1,E_2)},
$$
where the supremum is taken over all Borel spectral measures $E_1$ and $E_2$ on $\T$.

It is also well known (see \cite{BS4} and survey \cite{AP}) that under these assumptions, the following formula holds:
\bay
\label{DKBS}
f(U)-f(V)=\iint\limits_{\T\times\T}(\dg f)(\z,\t)\,dE_U(\z)(U-V)\,dE_V(\T),
\ey
where $E_U$ and $E_V$ are the spectral measures of $U$ and $V$.

Suppose that $E_1$ and $E_2$ are Borel spectral measures on locally compact topological spaces
$\X$ and $\Y$, at least one of which is separable and suppose that $\supp E_1=\X$ and $\supp E_2=\Y$.
Then taking into account Theorem 2.1 of  \cite{AP2}, we obtain from Theorem 2.2.4 of \cite{AP}
the following fact:

{\it Let $\Phi$ be a function on $\X\times\Y$ that is continuous in each variable.
Then \lb$\Phi\in\fM(E_1,E_2)$ if and only if it belongs to the Haagerup tensor product 
$C_{\rm b}(\X)\!\otimes_{\rm h}\!C_{\rm b}(\Y)$ of the spaces
$C_{\rm b}(\X)$ and $C_{\rm b}(\Y)$ of bounded continuous functions on
$\X$ and $\Y$, i.e., $\Phi$ admits a representation of the form
{\em\rf{predstava}, where
$\f_n\in C_{\rm b}(\X)$, $\psi_n\in C_{\rm b}(\Y)$ and the following inequalities hold:}
$$
\sum_n|\f_n|^2\le\|\Phi\|_{\fM(E_1,E_2)}
\quad\mbox{and}\quad\sum_n|\psi_n|^2\le\|\Phi\|_{\fM(E_1,E_2)}.
$$}

We proceed now to a general trace formula for double operator integrals.

Let $T$ be a trace class operator on Hilbert space, let
$E$ be a spectral measure on a $\s$-algebra of subsets of $\X$ and $\Phi\in\fM(E,E)$.
Let us compute the trace of the double operator integral
$$
\iint\Phi(x,y)\,dE(x)T\,dE(y).
$$
In \cite{BS4} the following formula was given:
\bay
\label{BStf}
\trace\left(\iint\Phi(x,y)\,dE(x)T\,dE(y)\right)=\int\Phi(x,x)\,d\mu(x),
\ey
where $\mu$ is the complex measure on the same $\s$-algebra that is defined 
by the equality
$$
\mu(\D)=\trace\big(TE(\D)\big).
$$

To justify the right-hand side of \rf{BStf} we should understand how we can
interpret the values of  the function $\Phi$ on the diagonal $\{(x,x):~x\in\X\}$.
In \cite{Pe+} the following interpretation of \rf{BStf} was given. We can define the trace ${\mathscr T}\Phi$
of a function $\Phi$ of  $\fM(E,E)$ on the diagonal by the equality
$$
({\mathscr T}\Phi)(x)\df\sum_n\f_n(x)\psi_n(x),
$$
where $\f_n$ and $\psi_n$ are functions in the representation \rf{predstava} of $\Phi$ in in terms of a 
Haagerup tensor expansion. Then the trace
${\mathscr T}\Phi$ of a function $\Phi$ of class $\fM(E,E)$ on the diagonal belongs to $L^\be(E)$
and does not depend on a representation \rf{predstava}. By  $\Phi(x,x)$ in \rf{BStf} we mean $({\mathscr T}\Phi)(x)$, see \cite{Pe5}, \S\:1.1.

Finally, suppose that $E$ is a Borel spectral measure on a locally compact topological space $\X$
and $\Phi$ is a function on $\X\times\X$ continuous in each variable.
Then the following assertion holds (see \cite{Pe5}):

\begin{thm}
\label{sleddvoi}
Let $E$ be a Borel spectral measure on a locally compact space $\X$
and let $\Phi$ be a function in $\fM(E,E)$. If $\Phi$ is continuous in each variable, then
{\em\rf{BStf}} holds for every operator $T$ in the trace class.
\end{thm}

In fact, it suffices to consider the case when $\supp E=\X$ and to consider a representation 
\rf{predstava} of $\Phi$ as an element of the Haagerup tensor product $C_{\rm b}(\X)\!\otimes_{\rm h}\!C_{\rm b}(\X)$.
It is easy to see that in this case $({\mathscr T}\Phi)(x)=\Phi(x,x)$, $x\in\X$.

\

\section{\bf Operator differentiability in the strong operator topology}
\setcounter{equation}{0}
\label{Difsil}

\

In this section for an operator Lipschitz function $f$ on $\T$,
we consider a problem of differentiability of the operator function
$$
t\mapsto f\big(e^{{\rm i}tA}U\big)
$$
in the strong operator topology,
where $U$ is a unitary operator and $A$ is a bounded
self-adjoint operator. Note that an analogue of the following theorem for functions of self-adjoint operators
was obtained in \cite{AP}, Th. 3.5.5; see also \S\:\ref{Samosop} of the present paper.

\begin{thm} 
\label{sildiftorn}
Let $f$ be an operator Lipschitz function on $\T$, let $U$ be a unitary operator and let $A$ be 
a bounded self-adjoint operator. Then

\bay
\label{vsilopto}
\lim_{t\to0}\frac1t\Big(f\big(e^{{\rm i}tA}U\big)-f(U)\Big)={\rm i}\int_\T\int_\T\t(\dg f)(\z,\t)\,dE_U(\z)A\,dE_U(\t),
\ey
where the limit is taken in the strong operator topology.
\end{thm}

Note that in  \cite{Pe2} formula  \rf{vsilopto} was obtained for functions $f$ in the Besov class $B_{\be,1}^1(\T)$,
in which case the limit on the left-hand side of \rf{vsilopto} exists in the operator norm.

We need the following auxiliary statement.

\begin{lem} 
\label{Xnun-n}
Let $\{X_n\}_{n\ge0}$ be a sequence in the space $\mB(\h)$ of bounded linear  operators
on a Hilbert space  $\h$ and let $\{u_n\}_{n\ge0}$ be a sequence in $\h$. Assume that
$$
\sum_{n\ge0}X_nX_n^*\le a^2I\quad\mbox{and}\quad
 \sum_{n\ge0}\|u_n\|^2\le b^2
$$
for nonnegative numbers $a$ and $b$. Then the series $\sum_{n\ge0}X_nu_n$
converges weakly and
$$
\Big\|\sum_{n\ge0}X_nu_n\Big\|\le ab.
$$
\end{lem}

\Pf Let $v\in\h$ and $\|v\|=1$. Then
$$
\sum_{n\ge0}|(X_nu_n,v)|=\sum_{n\ge0}|(u_n,X_n^*v)|\le
\Big(\sum_{n\ge0}\|u_n\|^2\Big)^{1/2}\Big(\sum_{n\ge0}\|X_n^*v\|^2\Big)^{1/2}\le ab,
$$
whence the desired result follows. $\bl$

\medskip

{\bf Proof of Theorem \ref{sildiftorn}.} 
As we have mentioned in \S\:\ref{Dvoiopi}, $f$ is differentiable on $\T$, 
the divided difference $\dg f$ is a Schur multiplier with respect to arbitrary Borel spectral measures 
$E_1$ and $E_2$, and
$$
\|\dg f\|_{\fM(E_1,E_2)}\le\|f\|_{\OL}.
$$

We have also observed in \S\:\ref{Dvoiopi} that there are sequences $\{\f_n\}_{n\ge0}$ and $\{\psi_n\}_{n\ge0}$ of continuous functions 
on $\T$ such that

{\rm a)} $\sum\limits_{n\ge0}|\f_n|^2\le\|f\|_{\OL(\T)}$ everywhere on $\T$,

{\rm b)} $\sum\limits_{n\ge0}|\psi_n|^2\le\|f\|_{\OL(\T)}$ everywhere on $\T$,

{\rm c)} $(\dg f)(\z,\t)=\sum\limits_{n\ge0}\f_n(\z)\psi_n(\t)$ for all  $\z$ and $\t$ from $\T$.

In view of identities \rf{doitproh} and \rf{DKBS}, we have to show that
$$
\lim_{t\to0}\frac1t\sum_{n\ge0}\f_n(e^{{\rm i}tA}U)(e^{{\rm i}tA}-I)U\psi_n(U)={\rm i}\sum_{n\ge0}\f_n(U)AU\psi_n(U).
$$
Here summation is considered in the weak topology while the limit is taken 
in the norm. Note that $\lim\limits_{t\to0}t^{-1}(e^{{\rm i}tA}-I)={\rm i}A$ 
in the operator norm. Thus, it suffices to prove that
$$
\lim_{t\to0}\sum_{n\ge0}\f_n(e^{{\rm i}tA}U)AU\psi_n(U)=\sum_{n\ge0}\f_n(U)AU\psi_n(U)
$$
in the strong operator topology.
In other words, we have to show that for every $u\in\h$,
$$
\lim_{t\to0}\sum_{n\ge0}(\f_n(e^{{\rm i}tA}U)-\f_n(U))AU\psi_n(U)u=\0,
$$
where summation is considered in the weak topology on $\h$,
and the limit is taken in the norm of $\h$. We assume that $\|u\|=1$ and $\|f\|_{\OL(\T)}=1$.
Then $\sum_{n\ge0}|\f_n|^2\le1$ and
$\sum_{n\ge0}|\psi_n|^2\le1$ everywhere on $\T$.

Put $u_n\df AU\psi_n(U)u$. We have
$$
\sum_{n\ge0}\|u_n\|^2\le\|A\|^2\sum_{n\ge0}\|\psi_n(U)u\|^2=
\|A\|^2\sum_{n\ge0}(|\psi_n|^2(U)u,u)\le\|A\|^2<+\be.
$$
Let $\e>0$. We can choose a positive integer $N$  such that 
$\sum_{n>N}\|u_n\|^2<\e^2$. Then it follows from Lemma  \ref{Xnun-n}
that
$$
\Big\|\sum_{n>N}(\f_n(e^{{\rm i}tA}U)-\f_n(U))u_n\Big\|\le2\e
$$
for all $t\in\R$. 

It is easy to see and it is well known that if  $h$  is a continuos function on $\T$, then the mapping
$$
U\mapsto h(U)
$$
is continuous in the operator norm on the set of unitary operators 
(it suffices to approximate $h$ by trigonometric polynomials).

Then
$$
\left\|\sum_{n=0}^N(\f_n(e^{{\rm i}tA}U)-\f_n(U))u_n\right\|
\le\|A\|\sum_{n=0}^N\Big\|\f_n\big(e^{{\rm i}tA}U\big)-\f_n(U)\Big\|<\e
$$
for all $t$ sufficiently close to zero. Thus,
$$
\Big\|\sum_{n\ge0}\big(\f_n(e^{{\rm i}tA}U)-\f_n(U)\big)u_n\Big\|<3\e
$$
for all $t$ sufficiently close to zero.  $\bl$

\

\section{\bf The trace formula and operator Lipschitzness}
\setcounter{equation}{0}
\label{Unitar}

\

In this section we establish the main result of this paper, which is the following theorem:

\begin{thm}
\label{osnrezu}
Trace formula {\em\rf{LKun}} holds for every operator Lipschitz function $f$ on $\T$ and for every
pair $(U,V)$ of unitary operators with trace class difference $U-V$.
\end{thm}

In the proof we use an idea by Birman and Solomyak which
was used in \cite{BS3} for their approach to a construction of the spectral shift function. However, in their paper stronger assumptions on $f$ were imposed.

\medskip

{\bf Proof of Theorem \ref{osnrezu}.} 
First of all, it is easy to see that under the condition $U-V\in\bS_1$, there exists a trace class self-adjoint operator $A$ such that $V=e^{{\rm i}A}U$. By Theorem \ref{sildiftorn}, the function $t\mapsto f\big(e^{{\rm i}tA}U\big)$
is differentiable in the strong operator topology and
\bay
\label{formdlyapro}
Q_s\df\frac{d}{dt}f\big(e^{{\rm i}tA}U\big)\Big|_{t=s}=
{\rm i}\int_\T\int_\T\t(\dg f)(\z,\t)\,dE_s(\z)A\,dE_s(\t),
\ey
where $E_s$ is the spectral measure of the unitary operator $V_s\df e^{{\rm i}sA}U$.

As we have observed in \S\:\ref{Dvoiopi}, the divided difference $\dg f$ is a Schur multiplier, hence by \rf{yadotse}
we have
$$
Q_s\in\bS_1\quad\mbox{and}\quad\sup_{s\in[0,1]}\|Q_s\|_{\bS_1}<\be.
$$

It follows from the definition of the function $s\mapsto Q_s$  that the function $s\mapsto Q_su$
is measurable for every $u\in\h$. Hence, the scalar function $s\mapsto (Q_su,v)$ is measurable
for every $u$ and $v$ in $\h$. This implies easily the measurability of  the function
$s\mapsto \trace(Q_sT)$ for every $T\in\mB(\h)$. Thus, the $\bS_1$-valued function $s\mapsto Q_s$
is weakly measurable, and so it is strongly measurable because the space $\bS_1$ is separable,
see, e.g., \cite{I}, Chapter V, \S\:4.

Now it follows from \rf{formdlyapro} that
$$
f(V)-f(U)=\int_0^1Q_s\,ds,
$$
where the integral is understood as the Bochner integral in the space $\bS_1$.

Then
$$
\trace\big(f(V)-f(U)\big)=\int_0^1\trace Q_s\,ds.
$$ 

By Theorem \ref{sleddvoi},
$$
\trace Q_s=\int_\T \z f'(\z)\,d\nu_s(\z),
$$
where $\nu_s$ is the complex Borel measure on $\T$ defined by 
$$
\nu_s(\D)\df\trace(E_s(\D)A)
$$
for Borel subsets $\D$ of $\T$.

We can identify the space ${\mathcal M}(\T)$ of the Borel complex measures on $\T$
with the space dual to the space  $C(\T)$ of continuous functions on $\T$.
Let us show that the function $s\mapsto\nu_s$ is continuous in the weak-$*$ topology of 
${\mathcal M}(\T)$.
Indeed, if  $h\in C(\T)$, then
$$
\int_\T h\,d\nu_s=\trace(h(V_s)A).
$$
As we have already observed in the proof of Theorem  \ref{sildiftorn},  the function $s\mapsto h(V_s)$
is continuous in the operator norm. This implies that the function $s\mapsto\nu_s$
is weakly continuous.

We define now the complex measure $\nu$ by the equality
$$
\nu=-\int_0^1\nu_s\,ds.
$$
Here the integral is understood as the integral of a continuous function in the weak-$*$ topology of ${\mathcal M}(\T)$.

Then
$$
\trace\big(f(U)-f(V)\big)=\int_\T\z f'(\z)\,d\nu(\z).
$$
On the other hand, we have
$$
\trace\big(f(U)-f(V)\big)=\int_\T f'(\z)\bs{\xi(\z)}\,d\z
$$
for trigonometrical polynomials $f$.
This implies that there exists a constant $c$  such that
$$
\z\,d\nu(\z)=\bs{\xi}(\z)\,d\z+c\,\z^{-1}\,d\z,
$$
which completes the proof of the theorem $\bl$

\medskip

Theorem  \ref{osnrezu} allows us to get the following amusing fact.

\begin{thm}
\label{podkrutka}
Let $f$ be an operator Lipschitz function on $\T$, and let $U$ and $V$ be unitary operators such that $U-V\in\bS_1$.
Then the function 
$$
\z\mapsto\trace\big(f(\z U)-f(\z V)\big),\quad\z\in\T,
$$
is continuous on $\T$.
\end{thm}

\Pf Let $f\in\OL(\T)$. Then for $\z\in\T$, we put  $f_\z(\t)=f(\z\t)$, $\t\in\T$.
We have
$$
\trace\big(f(\z U)-f(\z V)\big)=\trace\big(f_\z(U)-f_\z(V)\big)
=\int_\T f_\z'(\t)\bs{\xi}(\t)\,d\t,
$$
where $\bs{\xi}$ is the spectral shift function for the pair $(U,V)$. It remains to observe that the function
$$
\z\mapsto\int_\T f_\z'(\t)\bs{\xi}(\t)\,d\t,\quad\z\in\T,
$$
is continuous on $\T$ because  $\bs{\xi}$ is integrable and  $f'$ belongs to $L^\be$. $\bl$

\medskip

Note also that if one could find an independent proof of Theorem \ref{podkrutka},
then it would be easy to deduce from it Theorem \ref{osnrezu}.

\

\section{\bf An alternative approach to the case of self-adjoint operators}
\setcounter{equation}{0}
\label{Samosop}

\

The following result was used in  \cite{Pe5} to prove the Lifshits--Krein trace formula for operator Lipschitz functions
of self-adjoint operators:
let $A$ and $K$ be self-adjoint operators such that $K\in\bS_2$ and let
$f$ be a differentiable function on $\R$ with bounded derivative, then the function
$t\mapsto f(A+tK)-f(A)$ is differentiable in the $\bS_2$-norm and 
\bay
\label{kommunisticheskaya}
\frac d{dt}\big(f(A+tK)-f(A)\big)\Big|_{t=0}=
\int_\R\int_\R(\dg f)(x,y)\,dE_A(x)K\,dE_A(y).
\ey

An analogue of this statement for functions of unitary operators could be the following:
if $f$ is a differentiable function on $\T$ with bounded derivative, $U$ is a unitary operator and $A$ is 
a self-adjoint operator of class $\bS_2$, then the function $t\mapsto f\big(e^{{\rm i}tA}\big)U$
is differentiable in the $\bS_2$-norm and  \rf{vsilopto} holds.
Unfortunately, we do not know whether this is true.

Instead, we have used in this paper the differentiability of this function in the strong
operator topology in the case when $f$ is an operator Lipschitz function on $\T$, see Theorem \ref{sildiftorn}.

In \cite{AP}, Th. 3.5.6, the following analogue of Theorem  \ref{sildiftorn}  for functions of self-adjoint operators
was obtained:

{\it  Let $f$ be an operator Lipschitz function on  $\R$ and let $A$ and $K$ be self-adjoint operators such that $K$ is bounded. Then the function  $t\mapsto\big(f(A+tK)-f(A)\big)$ is differentiable
in the strong operator topology and} \rf{kommunisticheskaya} holds.

\medskip

This theorem allows us to obtain a new proof of the Lifshits--Krein trace formula for operator Lipschitz
functions of self-adjoint operators that does not use the result of \cite{KPSS} mentioned  above on the differentiability of operator functions in the Hilbert--Schmidt norm.

\

\

\noindent
\begin{tabular}{p{8cm}p{8cm}}
A.B. Aleksandrov  &  V.V. Peller \\
St.Petersburg Branch  & Department of Mathematics  \\
Steklov Institute of Mathematics  & Michigan State University\\
Fontanka 27   & East Lansing, Michigan 48824 \\
 191023 St-Petersburg  & USA\\
 Russia
\end{tabular}

\end{document}